\documentclass[12pt]{article}
\usepackage{mathrsfs}
\usepackage{algorithm}
\usepackage{algpseudocode}
\usepackage{amsmath}
\usepackage{amssymb}
\usepackage{bm}
\usepackage{booktabs}
\usepackage{color}
\usepackage{cite}
\usepackage{epsfig}
\usepackage{epstopdf}
\usepackage{flafter}
\usepackage{float}
\usepackage{ifthen}
\usepackage{multirow}
\usepackage{amsfonts}
\usepackage{makecell} 
\usepackage{pifont}
\usepackage{subfigure}
\usepackage{times}
\usepackage{threeparttable}
\usepackage{url} 
\usepackage{subfigure,graphicx}
\usepackage{float}
\usepackage{color,multirow}
\usepackage{makecell}
\usepackage{setspace}
\usepackage{cite} 
\usepackage{url} 
\usepackage{ifthen} 
\usepackage[T1]{fontenc}
\usepackage{graphicx}
\usepackage{changepage}
\usepackage{graphics}
\usepackage{booktabs}

\usepackage{amsfonts}
\usepackage{color}
\usepackage{amssymb}
\usepackage{amsmath}
\usepackage{indentfirst}

\def\sqr#1#2{{\vcenter{\vbox{\hrule height.#2pt
		\hbox{\vrule width.#2pt height#1pt \kern#1pt \vrule width.#2pt}
		\hrule height.#2pt}}}}
\def\signed #1{{\unskip\nobreak\hfil\penalty50
\hskip2em\hbox{}\nobreak\hfil#1
\parfillskip=0pt \finalhyphendemerits=0 \par}}
\def\endpf{\signed {$\sqr69$}}

\def\dbR{{\mathop{\rm l\negthinspace R}}}

\def\3n{\negthinspace \negthinspace \negthinspace }
\def\2n{\negthinspace \negthinspace }
\def\1n{\negthinspace }

\def\ds{\displaystyle}

\def\dbE{{\mathbb{E}}}
\def\dbF{{\mathbb{F}}}

\def\dbP{{\mathbb{P}}}

\def\dbR{{\mathbb{R}}}

\def\={\buildrel \triangle \over =}

\def\b{\beta}

\def\e{\varepsilon}

\def\l{\lambda}

\def\n{\nabla}

\def\t{\times}
\def\f{\varphi}
\def\th{\theta}
\def\o{\omega}

\def\i{\infty}
\def\ns{\noalign{\ss} }

\def\G{\Gamma}
\def\D{\Delta}

\def\O{\Omega}

\def\cC{{\cal C}}

\def\cF{{\cal F}}

\def\cW{{\cal W}}

\def\cl{{\cal l}}

\def\ss{\smallskip}
\def\ms{\medskip}

\def\q{\quad}
\def\qq{\qquad}

\def\esssup{\mathop{\rm esssup}}
\def\max{\mathop{\rm max}}

\def\inf{\mathop{\rm inf}}

\def\pa{\partial}

\def\cd{\cdot}

\def\div{\hbox{\rm div$\,$}}

\def\esssup{\hbox{\rm ess$\,$\rm sup$\,$}}

\def\as{\hbox{\rm a.s.{ }}}

\def\cl{\overline}

\def\deq{\mathop{\buildrel\D\over=}}

\def\|{\Big |}
\def\({\Big (}
\def\){\Big )}
\def\[{\Big[}
\def\]{\Big]}
\def\bel{\begin{equation}\label}
\def\ee{\end{equation}}
\def\bt{\begin{theorem}}
\def\bcd{\begin{condition}}
\def\ecd{\end{condition}}
\def\et{\end{theorem}}
\def\bc{\begin{corollary}}
\def\ec{\end{corollary}}
\def\bde{\begin{definition}}
\def\ede{\end{definition}}
\def\bl{\begin{lemma}}
\def\el{\end{lemma}}
\def\bp{\begin{proposition}}
\def\ep{\end{proposition}}
\def\br{\begin{remark}}
\def\er{\end{remark}}
\def\ba{\begin{array}}
\def\ea{\end{array}}
\def\ed{\end{document}}
\def\ns{\noalign{\ms}}
\def\ds{\displaystyle}

\def\square#1{\vbox{\hrule\hbox{\vrule height#1%
\kern#1\vrule}\hrule}}
\def\rectangle#1#2{\vbox{\hrule\hbox{\vrule height#1%
\kern#2\vrule}\hrule}}

\font\tenbb=msbm10 \font\sevenbb=msbm7 \font\fivebb=msbm5

\newfam\bbfam
\scriptscriptfont\bbfam=\fivebb \textfont\bbfam=\tenbb
\scriptfont\bbfam=\sevenbb

\newtheorem{lemma}{Lemma}[section]
\newtheorem{remark}{Remark}[section]

\newtheorem{theorem}{Theorem}[section]
\newtheorem{corollary}{Corollary}[section]

\newtheorem{definition}{Definition}[section]
\newtheorem{proposition}{Proposition}[section]
\newtheorem{condition}{Condition}[section]
\newtheorem{assumption}{Assumption}[section]

\makeatletter

\@addtoreset{equation}{section}
\makeatother

\allowdisplaybreaks

\begin{document}
\title{On Inverse Problems for Mean Field Games  with Common Noise  via Carleman estimate\thanks{This work is partially supported by the NSF of China under grants 12025105 and  12401586. }}
\author{  Zhonghua Liao\thanks{School of Mathematics,  Sichuan University,
Chengdu 610064,  China. E-mail address: zhonghualiao@yeah.net} \  and\  Qi L\"u\thanks{School
of Mathematics, Sichuan University, Chengdu
610064, China. E-mail address: lu@scu.edu.cn.}}
\date{}
\maketitle
\begin{abstract}
In this paper, we study two kinds of inverse problems for Mean Field Games (MFGs) with common noise. Our focus is on MFGs described by a coupled system of stochastic Hamilton-Jacobi-Bellman and Fokker-Planck equations.  Firstly, we establish the Lipschitz and H\"older stability for determining the solutions of a coupled system of stochastic Hamilton-Jacobi-Bellman and Fokker-Planck equations based on terminal observation of the density function. Secondly, we derive a uniqueness theorem for an inverse source problem related to the system under consideration.  The main tools to establish those results are two new Carleman estimates. 
\end{abstract}

{\bf Keywords: }   Mean field games, Carleman estimate, Lipschitz stability, H\"older stability, uniqueness

\section{Introduction}

\par  The theory of mean field games (MFGs for short) was first introduced by Lasry-Lions (\cite{JLPL}), as well as independently by Huang-Caines-Malham\'e (\cite{HCM}). It provides a tractable approximation to Nash equilibria in games involving a large number of players, where each individual's influence on the system is infinitesimally small but collectively significant. This theory has found applications in various fields, such as macroeconomics, crowd motions, finance,  and power grid models, where understanding the collective behavior of large groups is crucial (e.g.,\cite{RCFD1,HCM,JLPL,VNK}). Due to its wide-ranging applications,  MFGs theory  is extensively studied in recent years (see \cite{ABJFPY,Cardaliaguet2023,RCFD1, RCFD2, PCFDJLPL} and the rich references therein).

\par In this paper, we aim to investigate two inverse problems for mean-field games (MFGs) with common noise. 
To begin with, we  introduce some basic settings and notations for our model. 

Let $ (\O, \cF, \mathbb{F}, \dbP) $ with $\dbF=\{\cF_t\}_{t=0}^{+\i} $  be a complete filtered probability space on which a  $ n $-dimensional standard Brownian motion  $W(\cd)$ is defined and $\mathbb{F}$ is the natural filtration generated by  $ W(\cd) $, argumented by all $\dbP$-null sets.  

\par Let $\mathbf{H}$
be a Banach space and $T>0$.
\par $ \bullet  $  If 
$ f:\O\mapsto \mathbf{H} $ is Bochner integrable (w.r.t. $ \dbP $), we denote by $\dbE f$ the mathematical expectation of $f$, i.e.,
$$
\dbE f=\int_\O fd\dbP.
$$ 

$ \bullet  $ By $ L_{\cF_t}^p(\O; \mathbf H) $ ($t\in [0,T],~ p\ge 1 $) we denote all measurable random variables $ f $ such that $ \dbE \Vert f\Vert_{\mathbf H}^p<+\i $.

$ \bullet  $ By $ L_{\dbF}^\i(\O; C([0,T];\bm H))  $ we denote the space consisting of all $ \bm H $-valued, continuous, $ \dbF $-adapted processes $ X(\cd) $ such that $ \esssup_{\o\in \O} \Vert X(\cd)\Vert_{C([0,T]; \bm H)} <+\i $.

$ \bullet $ By $L^p_{\mathbb{F}}(0, T; {\bf H})$ ($ p\ge 1 $)
we denote  the  space consisting of all ${\bf H}$-valued $\mathbb{F}$-adapted  processes $X(\cdot)$ such that $\dbE
\Vert X(\cdot)\Vert_{L^p(0, T; {\bf H})}^p < +\infty$.

$ \bullet $ By $L^\infty_{\mathbb{F}}(0, T; {\bf H})$ we denote the
space of all ${\bf H}$-valued
$\mathbb{F}$-adapted bounded processes.


All of the aforementioned spaces are Banach spaces equipped with the canonical norm.

\par Let  $ \b\in [0,1] $  and $ \hat \b\= \frac{1+\b^2}{2} $. We consider the following coupled forward-backward stochastic parabolic equations: 
\begin{equation}\label{3}
\left\{
\ba{ll}
\ds d\rho -\hat \b \D \rho dt=\n\cd (\rho \n_p { H}(\cd , \cd, \cd, \n u;\rho))dt -\b \n \rho \cd dW(t)  &\mbox{ in } (0,T)\t \dbR^n,\\
\ns\ds du+\hat \b \D udt =-\b \div Udt- H(\cd,\cd ,\cd, \n u; \rho)dt+U\cd dW(t) & \mbox{ in } (0,T)\t \dbR^n,\\
\ns\ds u(T, \cd)=h(\cd),\q \rho(0,\cd)=\rho_0(\cd) & \mbox{ in } \dbR^n.
\ea 
\right.	
\end{equation}
Here $h\in L_{\cF_T}^2(\O; L^2(\dbR^n)),~\rho_0\in L^2 (\dbR^n)$, and $ H $ denotes the  Hamiltonian  as follows: for interaction functions $ \bm B(\cd,\cd,\cd,\cd): [0,T]\t \O \t \dbR^{2n}\to \dbR^n $, $b(\cd,\cd,\cd,\cd):[0,T]\t \O \t \dbR^{2n}\to \dbR$, $G(\cd,\cd, \cd): [0,T]\t \O \t \dbR^n \to \dbR$, $ F( \cd,\cd): \dbR^2\to \dbR $  and $ K(\cd, \cd, \cd,\cd): [0,T]\t\O\t  \dbR^{2n}\to \dbR $, 
\begin{equation}\label{a1-11}
\ba{ll}
\ds H(t, \o,  x, p; \rho)\!\!\!&\ds \= \inf_{\bm a\in \dbR^n}\Big\{\bm B(t,\o,x,\bm a)\cd p+ b(t,\o,x,\bm a)+G(t,\o,x)\\
\ns&\ds\qq\q  +F\(\int_{\dbR^n} K(t,\o, x,y)\rho(t,\o,y)dy, \rho(t,\o,x)\) \Big\}\\
\ns & \ds \= { B}(t,\o,x,p)+G(t,\o,x)+F\(\int_{\dbR^n} K(t,\o, x,y)\rho(t,\o,y)dy, \rho(t,\o,x)\),\\
\ns&\ds \qq\qq\qq\qq\qq\qq \qq\qq  (t,\o,x,p)\in [0,T]\t \O \t \dbR^{2n}, 
\ea 
\end{equation}
with $ { B}(t,\o, x,p)\= \inf\limits_{\bm a\in \dbR^n} \{\bm B(t,\o,x,\bm a)\cd p+b(t,\o,x,\bm a)\} $. 

In what follows, for the sake of brevity and clarity, we will omit the variable $ \o $ unless its inclusion is necessary to avoid ambiguity.

Next, we introduce the following assumptions for the functions appeared in \eqref{3} and \eqref{a1-11}:
\begin{assumption}\label{ass1}
Let $ F\in C^1(\dbR ^2;\dbR)$ and  $ { B} \in L_{\dbF}^\i (\O; C([0,T]; C^3(\dbR^{2n}; \dbR)))$. Moreover, For any $ M_1>1 $,  there exists a constant $ \cC(M_1)>0 $ such that for any $ (t,\o,x)\in [0,T]\t \O \t \dbR^n $, $ |p|\le M_1 $,  and for $ j',j,k=1,\cdots, n $,
\begin{equation}\label{b2-2}
|B_{p_jx_j}(t,x,p)|+|B_{p_jp_k}(t,x,p)|+|B_{p_jp_kx_k}(t,x,p)|+|B_{p_{j'}p_jp_k}(t,x,p)| \le \cC(M_1).
\end{equation}
\end{assumption}
\begin{assumption}\label{ass2}
Let $ K\in L_\dbF^\i(0,T; L^2(\dbR^{2n})) $ and  there exists a   positive constant $ M_2 $ such that 
$ \Vert K\Vert_{L_\dbF^\i(0,T; L^2(\dbR^{2n}))}\le M_2 $.
\end{assumption}  

The system (\ref{3}) describes a mean field models with common noise (see \cite[Section 4.1]{PCFDJLPL}). 
In the system \eqref{3}, the first equation corresponds to the Fokker-Planck (FP) equation for the density function, while the second equation represents the stochastic Hamilton-Jacobi-Bellman (HJB) equation for the value function. The solution to   \eqref{3} is given by a triple $(\rho, u, U)$, where $\rho$ denotes the density function, $u$ represents the value function, and $U$ stands for the correction terms in the backward stochastic parabolic equation.  

\begin{remark}
In practical applications, the term $\int_{\mathbb{R}^n} K(x,y) \rho(y)dy$ describes the interaction between an individual player and the surrounding environment. For example, in a traffic congestion model, each driver's strategy is influenced not only by their own position but also by the density of traffic in the vicinity. Similarly, in a bird migration scenario, a bird's migration strategy is affected by the presence of other birds nearby. 
A common scenario occurs when $K(t, x,y) = K(t, x-y)$ serves as a convolution operator kernel. In such cases, Assumption \ref{ass2} necessitates that $K(\cdot)$ has compact support and belongs to $L^2(\mathbb{R}^n)$.
\end{remark}

As one can see, the FP-HJB system  (\ref{3}) consists a coupled forward-backward  parabolic equations. It is natural to expect the well-posedness after given $ u(T,\cd) $ and $ \rho(0,\cd) $. However, one may fails because the uniqueness of this system is quite rare, and have only been proven under strong conditions, like, e.g., a monotonicity assumption \cite{JLPL}.  There are   examples showing that (\ref{3}) may have two classical solutions even for smooth Hamiltonian (\cite{MBMF}).  Consequently,  the first inverse problem we are interested in  is that  whether we  can weaken those assumptions, but restore the uniqueness  by adding an additional measurement. More precisely,  besides $ u(T,\cd) $ and $ \rho(0,\cd) $, we also assumed the availability of  $ \rho(T,\cd ) $, and our goal is to determine whether this extra data can restore  uniqueness. Furthermore, since the measurements might be given with a noise/error. Therefore, it is also important to get the corresponding stability estimate. In practical terms, this inverse problem represents a scenario where the games has already happened, and we are interested that whether the process  depends on the terminal state  continuously. Specifically, we study the following questions:

\ss

\par {\bf (IP1) } Can we determine the solution by proper measurements?  Further,  is the solution continuously dependent on  the measurements?

\ss

Next, we consider an inverse source problem for the FP-HJB system  (\ref{3}).
Let us assume that the source term in (\ref{3}) has the form  $ G(t,\o, x)=r(t,\o, x')R(t,x) $  with  an unknown  $ r \in L_\dbF^2(0,T; H^1(\dbR^{n-1})) $ which is independent of the first component of $ x\=(x_1, x')\in \dbR^n $ ( where $ x'=(x_2,\cdots,x_n)\in \dbR^{n-1} $). 

\ss

\par {\bf (IP2)} Let $ R(\cd,\cd)$ in $G(\cd, \cd,\cd)$ be given. Can we determine the source function $ r(t,x') $, $ (t,x')\in (0,T)\t \dbR^{n-1} $ by means of the observation of $ \rho(T,x) $ and  Lateral Cauchy data on  $ x_1=l_1 $ and $ x_1=l_2 $ for some $ l_1, l_2\in \dbR $? 

\ss

Problem \textbf{(IP2)} naturally arises in many Mean Field Games (MFGs). For instance, in numerous control problems in economic theory, we may have partial information about the cost function without knowing its exact form. In such scenarios, we aim to determine the unknown part through additional observations of the solution to the FP-HJB system (\ref{3}).

\ss

It is important to study the associated inverse problems in the field of MFGs. For example, in economics, solving proper inverse problems linked to market competition and resource allocation can enrich our comprehension and forecasting of market behaviour and outcomes. Similarly, in transportation, investigating MFG inverse problems can facilitate the optimization of traffic flow control strategies, thereby enhancing the efficiency and sustainability of road networks. In the realm of financial systems, exploring inverse problems in MFGs can provide insights into optimizing portfolio management strategies and risk assessment, leading to more robust investment decisions. 
Furthermore, in the field of environmental sustainability, exploring MFG inverse problems can assist in improving decision-making processes for resource management and conservation efforts, leading to a more balanced and eco-friendly approach to development.

Due to their significant applications, inverse problems for MFGs have attracted considerable attention in recent years (e.g., \cite{PKSASS, LDWLSO, MVKYA, MVK1, MVK2, MVKJLHL, MVKJLHL1, MVKJLZY, HLCMSZ, KRNSKW, KRNSKWHZ} and references therein). To the best of our knowledge, the pioneering work on the inverse problem {\bf (IP1)} for MFGs without common noise is presented in \cite{MVKYA}, where  the FP-HJB system are coupled forward-backward deterministic parabolic equations, rather than coupled forward-backward stochastic parabolic equations.  A key aspect of their method is the establishment of two innovative Carleman estimates for the FP-HJB system. These estimates lead to Lipschitz stability results. Subsequently, in \cite{MVKJLHL}, a similar Carleman estimate is utilized to derive a H\"older-type stability estimate with fewer measurements. Since then, various other Carleman estimates for  FP-HJB system  have been developed and applied to different inverse problems (e.g., \cite{MVK1, MVK2, MVKJLHL1, MDHLGZ}). However, it is important to note that all these works are focused exclusively on MFGs without common noise.

Common noise appears in many practical control systems, affecting all agents simultaneously (e.g.,\cite{RCFD2}). For example, in financial markets, the effect of a shared market factor or index can be represented as common noise influencing investors' decision-making; and in models of traffic flow, variables such as weather conditions or road constructions can be considered as common noise impacting motorists' driving behaviours. Consequently, incorporating common noise into MFGs not only increases complexity but also captures crucial phenomena. Hence, exploring the inverse problems of MFGs with common noise is desirable.


Borrow some ideas from \cite{MVK1, MVK2, MVKJLHL1,MVKYA,MVKJLHL,MDHLGZ}, to answer   {\bf Problem  (IP1)} and {\bf Problem  (IP2)}, we aim to  establish a suitable Carleman estimates for the FP-HJB system (\ref{3}) to solve the desired inverse problem. Due to their wide applications, Carleman estimates for stochastic parabolic equations have recently attracted a lot of attention (e.g., \cite{Liao2024,QL,QLXZ,XZQL,STXZ,Wu2020,Wu2022,GY} and the reference therein). Nevertheless, to establish proper Carleman estimates to solve {\bf Problem  (IP1)} and {\bf Problem  (IP2)}, we encounter  several challenges: 
\begin{itemize}
\item  Firstly, when considering the Carleman estimate for (\ref{3}),  unlike the deterministic FP-HJB system,  the diffusion term ``$ -\b \n\rho \cd dW(t) $" introduces additional undesirable  terms which have to be got rid of. 
\item Secondly,   there are non-homogeneous correction terms ``$ U\cd dW(t) $'' and ``$ -\b \div U\cd dt $" in the HJB equation, which add complexity to the analysis and require meticulous scrutiny.
\item Thirdly, solutions to stochastic parabolic equation may not be continuously differentiable with respect to the time variable. Hence, one cannot take the time derivative to reduce the non-homogeneous term in the stochastic equation under consideration to the initial datum of another equation as people did in the deterministic situation.
\end{itemize}
These difficulties have already been observed in \cite{STXZ} when the authors study Carleman estimate for stochastic parabolic equation (see \cite[Remarks 2.1 and 6.1]{STXZ}). For the HJB equation in \eqref{3}, the situation is even worse due to the appearance of the convolution interaction term $ \int_{\dbR^n} K(t, x, y)\rho(t, y)dy $. Fortunately, for our special inverse problems, after careful analysis and calculation,  we are able to address the gradient drift term in FP-HJB system in our Carleman estimates. This is primarily due to the relationship between the coefficients $ \b $ and $ \hat \b $ in  (\ref{3}).

Next, we introduce an  a prior bound of the solution $ (u, \rho) $ which is a common condition when investigating  stability estimate for  nonlinear equations. 
\begin{assumption}\label{ass3}
There exists a  positive constant $ M_3 $ such that
\begin{equation}\label{24}
\Vert u\Vert_{L_\dbF^\i(0,T; W^{2,\i}(\dbR^n))}+\Vert \rho\Vert_{L^\i_\dbF (0,T; W^{1, \i}(\dbR^n))}\le M_3. 
\end{equation}	
\end{assumption}

Now we  present  our main results in this paper, which provide answers to {\bf Problem  (IP1)} and {\bf Problem  (IP2)}.  For simplicity of notations, let us denote $ M\=\max\{ 1, \cC(M_3), M_2, M_3\} $. 

First, we  have the following Lipschitz stability estimate for {\bf Problem  (IP1)}.

\begin{theorem}\label{theorem1}
Let  Assumptions \ref{ass1}--\ref{ass3} hold. Suppose that  $ (\rho_i, u_i, , U_i)\in  L_\dbF^2(0,T, H^1 (\dbR^n))\t L_\dbF^2(0,T;$ $H^1(\dbR^n))\t L_\dbF^2(0,T; L^2(\dbR^n))$ ($ i=1,2 $) are two solutions of (\ref{3}).
Then there exists a constant $ C=C(T,M, F)>0 $, such that  for any  $ \b\in [0,1] $, we have
\begin{equation}\label{27}
\begin{aligned}
	& \Vert (u_1-u_2, \rho_1-\rho_2)\Vert_{L_\dbF^2(0,T, H^1(\dbR^n))\t L_\dbF^2(0,T; H^1(\dbR^n))}\\ & \le C\(\dbE \Vert u_1(T)-u_2(T)\Vert_{H^1(\dbR^n)}+\dbE \Vert \rho_1(T) - \rho_2(T) \Vert_{L^2 (\dbR^n)} + \Vert \rho_1(0) - \rho_2(0) \Vert_{H^1(\dbR^n)} \).
\end{aligned}
\end{equation}
\end{theorem}

In Theorem \ref{theorem1}, we employ three measurements   to determine the solution. Next, we only use two measurements $ u_i(T,\cd) $ and $ \rho_i(T,\cd) $. In this case, we can get the following H\"older type stability estimate for {\bf Problem  (IP1)}. 

\begin{theorem}\label{theorem2}
Let $ \e\in (0,T) $   and  Assumptions \ref{ass1}-\ref{ass2} hold. Suppose that  $ (\rho_1, u_1, U_1) $ and $ (\rho_2, u_2, U_2) $ are two solutions of (\ref{3}) satisfying Assumption \ref{ass3}. Then there exist  constants  $ \eta=\eta(T, M, F, \e) $ and $ C=C(T,M, F, \e) $, such that for any  $ \b\in [0,1] $, we have
\begin{equation}\label{a3}
\ba{ll}
\ds \Vert \rho_1-\rho_2\Vert_{L^2_\dbF (\e, T; H^1(\dbR^n))}+\Vert u_1-u_2\Vert_{L^2(\e, T; H^1(\dbR^n))}\\
\ns\ds \le 	C \Vert \rho_1(0)-\rho_2(0)\Vert_{H^1(\dbR^n)}^{1-\eta}\big[\Vert u_1(T)-u_2(T)\Vert_{L^2_{\cF_T}(\O; H^1(\dbR^n))}+\Vert \rho_1(T) - \rho_2(T)\Vert_{L_{\cF_T}^2(\O; L^2(\dbR^n))}\big]^\eta.
\ea 
\end{equation}
\end{theorem}

\par Next, we consider {\bf Problem  (IP2)}.
Without loss of generality, we set $ l_1=0 $ and $ l_2=1 $. Denote $ G\= (0,1)\t \dbR^n $. Suppose that $ (\rho_i, u_i, U_i)\in L_\dbF^2 (0,T; H^2(G))\t L_\dbF^2(0,T; H^1(G))\t L_\dbF^2(0,T;  L^2(G)) $ ($ i=1,2 $) are two solutions of system (\ref{3}) with  the same terminal cost $ u_1(T,\cd)=u_2(T, \cd) $, but two different source terms $  r_{1}, r_{2} $. We have the following uniqueness result for the above problem.  

\begin{theorem}\label{theorem3}
Suppose Assumptions \ref{ass1}-\ref{ass3} hold and $\Vert K_{x_1}\Vert_{L_\dbF^\i(0,T; L^2(G\t G))}\le M $. 
Suppose   $ R\in C^2([0,T]\t \cl G) $  and $ \inf_{x\in \cl G}|R(t,x)|> 0 $. If $ \rho_1(T, \cd)=\rho_2(T, \cd)$ on $ \dbR^n $, $\dbP$-{\rm a.s.} and  
\begin{equation}\label{d1}
\begin{cases}\ds
	\rho_1=\rho_2,\q  u_1=u_2,\q \rho_{1x_1}=\rho_{2x_1},  \\
	\ns\ds u_{1x_1}=u_{2x_1}, \q \rho_{1x_1x_1}=\rho_{2x_1x_1},\q u_{1x_1x_1}=u_{2x_1x_1},
\end{cases} \q \mbox{ on } [0,T]\t \pa G,\;\dbP\mbox{-}\as,
\end{equation}
then $ r_1(t,x')=r_2(t,x') $ for all $ (t,x')\in [0,T]\t \dbR^{n-1} $, $\dbP$-{\rm a.s.}
\end{theorem}
\begin{remark}
In practice, we can conduct polls of game players to get the boundary data. Typically, these polls are not limited to the boundary itself, but conducted in a small neighbourhood of the boundary. Hence, it is possible to approximately figure out all the data in (\ref{d1}). Notably,  if we have  prior knowledge that the source term $ r $ is independent of time variable, our result implies that, theoretically, one can perform an experiment over a brief time interval to reconstruct the source term.
\end{remark}

The rest of the paper is organized as follows. In Section 2,, we establish two Carleman estimates for  forward and backward stochastic parabolic operators, which play key roles in proving Theorems \ref{theorem1}--\ref{theorem3}. In Sections 3 and 4, we give  proofs of our main results.

\section{Two Carleman Estimates for FP-HJB system}

In this section, we establish two Carleman estimates for the system (\ref{3}), which are crucial in solving {\bf Problem  (IP1)} and {\bf Problem  (IP2)}. To begin with,   we introduce the following weight function:
\begin{equation}\label{4}
\th\= e^\ell, \q \ell\= \l \f,\q \f\=(t+2)^\mu,
\end{equation}
for some parameters $ \l, \mu > 0 $. 

Let $ w \in  L_\dbF^2(0,T; H^2(\dbR^n))$ be a solution to 
\begin{equation}\label{5*}
dw +\hat\b \D  wdt=(f_1-\b \div \bm f_2)dt+\bm f_2\cd dW(t),	
\end{equation}
where $ \hat\b$ and $\b $ are given in (\ref{3}), $ f_1\in L^2_\dbF(0,T, L^2(\dbR^n)) $, and  $ \bm f_2\in L_\dbF^2(0,T; H^1 (\dbR^n; \dbR^n))$.
We have the following Carleman estimate.

\begin{theorem}\label{th1}
For any $ \l, \mu>0 $ and  $ \b \in [0,1] $, it holds that
\begin{eqnarray}\label{11} 
&&   \frac{1}{2} \dbE\int_0^T\int_{\dbR^n}  \th^2 \[ \hat\b^2  |\D w|^2 +\l^2\mu^2 (t+2)^{2\mu-2}  w^2+2\l\mu(t+2)^{\mu-1} |\n w|^2 \] dxdt \nonumber \\
&&  + \dbE\int_0^T \int_{\dbR^n} \[ \l\mu(t+2)^{\mu-1} \th|\bm f_2|^2+  \frac{1-\b^2}{2} \th^2|\div \bm f_2|^2\]dxdt\\
&&  \le  2 \dbE \int_0^T \int_{\dbR^n} \th^2 f_1^2dxdt+ e^{2\l (T+2)^\mu } \dbE \int_{\dbR^n}  \big[ \l \mu (T+2)^{\mu-1} w(T,x)^2+\hat \b  |\n w(T,x)|^2 \big]dx.  \nonumber
\end{eqnarray}
\end{theorem}

\par {\bf Proof of Theorem \ref{th1}. } Set $ v=\th w $. Then
\begin{equation}\label{6}
\th (dw+\hat \b \D wdt)= dv -\l\mu(t+2)^{\mu-1} vdt+ \hat\b \D vdt.	
\end{equation} 
Multiplying (\ref{6}) by $ 2\hat \b \D v-2\l\mu (t+2)^{\mu-1} v $,  using It\^o's formula and (\ref{5*}),  we get that
\begin{eqnarray}\label{7}
&& 
\big[2\hat \b \D v-2\l\mu (t+2)^{\mu-1} v \big] \th  (dw+\hat \b \D wdt) \nonumber \\
&& =-d\big[\l\mu(t+2)^{\mu-1} v^2+\hat \b |\n v|^2 \big]+2\hat \b \n\cd (\n v dv)+\l\mu(t+2)^{\mu-1}(dv)^2+\hat \b|d\n v|^2	\nonumber\\
&& \q +\l\mu(\mu-1)(t+2)^{\mu-2} v^2dt+ 2\big[\hat \b \D v-\l\mu (t+2)^{\mu-1} v\big]^2dt.
\end{eqnarray}
Integrating both sides of (\ref{7}) over $(0,T)\times \dbR^n$,  noting that
\begin{equation}\label{a1-3} 
\ba{ll}
\ds \dbE\int_{\dbR^n}\int_0^T\Big[ \l\mu(t+2)^{\mu-1} (dv)^2+\hat \b |d \n v|^2\Big]dx\\
\ns\ds =\dbE \int_0^T \int_{\dbR^n} \l\mu(t+2)^{\mu-1}\th^2  |\bm f_2|^2+\hat\b \sum_{j=1}^n\th^2  |\bm f_{2x_j}|^2dxdt,
\ea 
\end{equation}
we get
\begin{eqnarray}\label{8}
&&\dbE \int_0^T\int_{\dbR^n} \Big[\l\mu(t+2)^{\mu-1} \th^2|\bm f_2|^2 + \hat \b \th^2 |\div  \bm f_2|^2 \Big]dxdt  \nonumber\\
&& +\dbE\int_0^T\int_{\dbR^n} \Big\{\l\mu(\mu-1)(t+2)^{\mu-2} v^2 + 2\big[\hat \b \D v-\l\mu (t+2)^{\mu-1} v\big]^2\Big\}dxdt \\
&& \le \dbE \int_0^T\int_{\dbR^n} \Big\{\frac{2}{3} \th^2  \big( 3f_1^2 +  \frac{3}{2}\b^2 |\div \bm f_2|^2 \big) +\frac{3}{2}\big[\hat \b \D v-\l\mu (t+2)^{\mu-1} v\big]^2\Big\} dxdt  \nonumber\\
&& \q + e^{2\l (T+2)^\mu } \dbE \int_{\dbR^n} \big[  \l \mu (T+2)^{\mu-1} w(T,x)^2+\hat \b  |\n w(T,x)|^2 \big]dx. \nonumber
\end{eqnarray}
Recalling that  $\b\in [0,1]$ and $ \hat \b=\frac{1+\b^2}{2} $,  we obtain from \eqref{8} that
\begin{equation}\label{9}
\ba{ll}
\ds \dbE\int_0^T \int_{\dbR^n} \[ \l\mu(t+2)^{\mu-1} \th|\bm f_2|^2+  \frac{1-\b^2}{2} \th^2|\div \bm f_2|^2\]dxdt\\
\ns\ds +\dbE\int_0^T\int_{\dbR^n} \Big\{\l\mu(\mu-1)(t+2)^{\mu-2} v^2 +   \frac{1}{2} \big[\hat \b \D v-\l\mu (t+2)^{\mu-1} v\big]^2\Big\}dxdt\\
\ns\ds \le  2 \dbE \int_0^T\int_{\dbR^n} \th^2  f_1^2  dxdt + e^{2\l (T+2)^\mu } \dbE \int_{\dbR^n} \big[  \l \mu (T+2)^{\mu-1} w(T,x)^2+\hat \b  |\n w(T,x)|^2 \big]dx.
\ea  
\end{equation}
Noting that
\begin{equation*}\label{10}
\ba{ll}
\ds \big[\hat \b \D v-\l\mu(t+2)^{\mu-1}v\big]^2 \\
\ns\ds =\hat \b^2|\D v|^2+\l^2\mu^2(t+2)^{2\mu-2} v^2-2\l\mu(t+2)^{\mu-1} \n\cd (v \n v)+2\l\mu(t+2)^{\mu-1}|\n v|^2,
\ea 
\end{equation*}
the inequality (\ref{11}) follows from (\ref{9}).
\endpf 

\ms

Let $ p\in L_\dbF^2 (0,T; H^2(\dbR^n)) $ be a  solution of 
\begin{equation}\label{12*}
dp -\hat \b\D p  dt=g_1dt- \b \n p\cd dW(t)
\end{equation}
for some $ g_1 \in L_\dbF^2 (0,T; L^2(\dbR^n))$. We have the following Carleman estimate.
\begin{theorem}\label{th2}
For any $ \l\ge 0 $, $  \mu\ge 144(T+2)^2 $ and  $ \b\in [0,1] $, we have
\begin{eqnarray}\label{20}
&&\frac{1}{4}\l\mu^2 \dbE \int_0^T\int_{\dbR^n}  (t+2)^{\mu-2}\th^2 p^2 dxdt+\sqrt \mu \dbE\int_0^T\int_{\dbR^n} \th^2 |\n p|^2 dxdt \nonumber\\
&& \le 2 \dbE \int_0^T \int_{\dbR^n} \th^2g_1^2dxdt+ 2 \l\mu(T+2)^{\mu-1}e^{2\l(T+2)^\mu} \dbE \int_{\dbR^n} p(T,x)^2dx\\
&&\q +\hat \b e^{2\l2^\mu}\dbE \int_{\dbR^n}|\n p(0,x)|^2dx+\sqrt \mu e^{2\l 2^{\mu}} \dbE \int_{\dbR^n} p(0,x)^2dx. \nonumber
\end{eqnarray}
\end{theorem}

{\bf Proof of Theorem \ref{th2}. } Set $ q= \th p $. Then 
\begin{equation}\label{13}
\th \big(dp-\hat \b \D pdt\big)=d q-\l\mu(t+2)^{\mu-1} qdt-\hat \b \D qdt.
\end{equation}
Multiplying (\ref{13}) by $- 2[(1+\b^2)\l\mu(t+2)^{\mu-1} q + \hat \b \D q ]$ and using It\^o's formula, we obtain that
\begin{eqnarray}\label{14} 
&& - 2\big[(1+\b^2)\l\mu(t+2)^{\mu-1} q + \hat \b \D q \big]\th \big(dp-\hat \b \D pdt\big)\nonumber\\
&&  =2\big[(1+\b^2)\l\mu(t+2)^{\mu-1} q + \hat \b \D q \big]^2dt \nonumber\\ && \q -2\big[(1+\b^2)\l\mu(t+2)^{\mu-1} q + \hat \b \D q \big]\big[dq+\b^2\l\mu(t+2)^{\mu-1} qdt\big]\\
&&  =2\big[(1+\b^2)\l\mu(t+2)^{\mu-1} q + \hat \b \D q \big]^2dt -d\big[(1+\b^2)\l\mu(t+2)^{\mu-1}q^2\big]\nonumber\\
&&  \q +(1+\b^2)\l\mu(\mu-1)(t+2)^{\mu-2}q^2dt+(1+\b^2)\l\mu(t+2)^{\mu-1}(dq)^2\nonumber\\
&&  \q -\n\cd (2\hat\b dq \n q)+\hat \b d(|\n q|^2 )-\hat \b |d \n q|^2-2\b^2 (1+\b^2)\l^2\mu^2(t+2)^{2\mu-2}q^2\nonumber\\
&&  \q -\b^2(1+\b^2)\l\mu(t+2)^{\mu-1}\D q qdt. \nonumber
\end{eqnarray}
Integrating both sides of (\ref{14}) over $(0,T)\times \dbR^n$,  noting that $ \hat \b=\frac{1+\b^2}{2} $ and $ \mu\ge 2 $,  by (\ref{12*}), we get that
\begin{equation*}\label{a1-4}
\ba{ll}
\ds \dbE\int_0^T \int_{\dbR^n}\Big[ (1+\b^2)\l\mu(t+2)^{\mu-1} (dq)^2-\hat \b |d\n q|^2\Big]dt\\
\ns\ds =\dbE \int_0^T\int_{\dbR^n}\Big[(1+\b^2)\l\mu(t+2)^{\mu-1} \b^2 |\n q|^2-\hat\b \b^2 \sum_{j,k=1}^n q_{x_jx_k}^2\Big]dxdt,
\ea 	
\end{equation*} 
which yields that
\begin{eqnarray}\label{15}
&&\dbE \int_0^T\int_{\dbR^n} \[ \frac{1+\b^2}{2}\l\mu^2(t+2)^{\mu-2}\th^2 p^2+\b^2(1+\b^2)\l\mu(t+2)^{\mu-1}|\n q |^2\]dxdt \nonumber \\
&& \q + \dbE \int_0^T\int_{\dbR^n}  2\big[(1+\b^2)\l\mu(t+2)^{\mu-1} q + \hat \b \D q \big]^2dxdt  \nonumber\\
&& \q +(1+\b^2)\l\mu 2^{\mu-1} e^{2\l 2^{\mu}}\dbE \int_{\dbR^n} p(0,x)^2dx+\hat \b e^{2\l(T+2)^{\mu}}\dbE \int_{\dbR^n}|\n p(T,x)|^2dxdt \nonumber\\
&& \q+\dbE \int_0^T\int_{\dbR^n}\big[ -2\b^2(1+\b^2)\l^2\mu^2(t+2)^{2\mu-1}q^2-\b^2(1+\b^2)\l\mu(t+2)^{\mu-1} \D q q \big] dxdt  \nonumber\\
&& \le -\dbE \int_0^T\int_{\dbR^n}  2\big[(1+\b^2)\l\mu(t+2)^{\mu-1} q + \hat \b \D q \big]\th \big(dp-\hat \b \D pdt\big)dx  \\
&& \q +\hat \b \b^2\dbE \int_0^{T}\int_{\dbR^n} \th^2  \sum_{j,k=1}^n p_{x_jx_k}^2dxdt \nonumber\\
&& \q +(1+\b^2)\l\mu(T+2)^{\mu-1} e^{2\l(T+2)^{\mu}} \dbE \int_{\dbR^n}  p(T,x)^2dx+\hat \b e^{2\l 2^{\mu}}\dbE \int_{\dbR^n} |\n p(0,x)|^2dx. \nonumber
\end{eqnarray}
Direct computation yields
\begin{equation}\label{34-2}
\ba{ll}
\ds \b^2(1+\b^2)\l\mu(t+2)^{\mu-1} |\n q |^2-2\b^2(1+\b^2)\l^2\mu^2(t+2)^{2\mu-1}q^2\\
\ns\ds \q -\b^2(1+\b^2)\l\mu(t+2)^{\mu-1} \D q q +\frac{2\b^2}{1+\b^2}\big[(1+\b^2)\l\mu(t+2)^{\mu-1} q + \hat \b \D q \big]^2\\
\ns\ds =\n\cd \big[\b^2(1+\b^2)\l\mu(t+2)^{\mu-1}  q \n q \big]+\hat \b \b^2 (\D q)^2,
\ea 
\end{equation}
and
\begin{equation}\label{34*}
\begin{aligned}
\sum_{j,k=1}^n \int_{\dbR^n}  p_{x_jx_k}^2 dx &= \sum_{j,k=1}^n\int_{\dbR^n}   \big[(p_{x_j} p_{x_jx_k})_{x_k}-(p_{x_j}p_{x_kx_k})_{x_j}+p_{x_jx_j}p_{x_kx_k}\big]dx \\ &=\int_{\dbR^n} (\D p)^2dx.	
\end{aligned}
\end{equation}
Recalling that $ \b\le 1 $, we have $ \frac{2\b^2}{1+\b^2}\le 1 $. By (\ref{12*}), (\ref{15}), (\ref{34-2}) and (\ref{34*}), we obtain that
\begin{eqnarray}\label{34-1} 
&& \!\!\!\!\dbE \int_0^T\int_{\dbR^n}  \frac{1}{2}\l\mu^2(t+2)^{\mu-2}\th^2 p^2 dxdt+\l\mu 2^{\mu-1} e^{2\l 2^{\mu}}\dbE \int_{\dbR^n} p(0,x)^2dx \nonumber\\
&&   \q +\hat \b e^{2\l(T+2)^{\mu}}\dbE \int_{\dbR^n}|\n p(T,x)|^2dx \\
&& \!\!\!\!\le \dbE\! \int_0^T\! \int_{\dbR^n}\! \th^2g_1^2dxdt \!+\! 2 \l\mu(T\!+\!2)^{\mu-1}e^{2\l(T+2)^\mu} \dbE\! \int_{\dbR^n}\! p(T,x)^2dx \nonumber \\
&& \q + \hat \b e^{2\l2^\mu}\dbE\! \int_{\dbR^n}\!|\n p(0,x)|^2dx. \nonumber
\end{eqnarray}
By It\^o's formula, we get that
\begin{equation}\label{16}
\ba{ll}
\ds 2 \th^2 p(dp-\hat \b \D pdt)\\
\ns\ds  = d( \th^2 p^2)-\th^2(dp)^2 -2 \l\mu(t+2)^{\mu-1} \th^2 p^2dt-\n\cd (2\hat \b \th^2 p \n p )dt+2\hat\b\th^2  |\n p|^2dt.
\ea 
\end{equation}
Integrating \eqref{16} over $ (0,T)\t \dbR^n $ and taking mathematical expectation, we obtain that
\begin{eqnarray}\label{19}
&& \sqrt \mu \dbE\int_0^T\int_{\dbR^n} \th^2 |\n p|^2 dxdt \nonumber\\
&& \le 3 \l \mu^{\frac{3}{2}}  \dbE \int_{0}^T \int_{\dbR^n}  (t+2)^{\mu-1}\th^2 p^2dxdt+\dbE\int_0^T\int_{\dbR^n} \th^2 g_1^2dxdt \\
&& \q +\sqrt \mu e^{2\l 2^{\mu}} \dbE \int_{\dbR^n} p(0,x)^2dx. \nonumber
\end{eqnarray}
Noticing that $ \sqrt \mu \ge 12(T+2) $, we have
\begin{eqnarray}\label{b1-1-1}
&&\[\frac{1}{2}\l\mu^2-3\l\mu^{\frac{3}{2}}(T+2)\]\dbE \int_0^T \int_{\dbR^n}  (t+2)^{\mu-2} \th^2 p^2 dxdt \nonumber\\
&&\ge\frac{1}{4}\l\mu^2 \dbE \int_0^T\int_{\dbR^n}  (t+2)^{\mu-2}\th^2 p^2 dxdt.
\end{eqnarray}
Combining (\ref{19}) and \eqref{b1-1-1} with (\ref{34-1}), we get (\ref{20}).
\endpf 

Next, we  introduce two Carleman estimates for solutions to stochastic parabolic equations on $(0,T)\t G$. 
\begin{theorem}\label{th3}
Let $ w \in  L_\dbF^2(0,T; H^2(G))$ be a  solution of 
\begin{equation}\label{5}
\begin{cases}\ds 
	dw +\hat\b \D  wdt=(f_1+\bm f_3\cd \bm f_2-\b \div \bm f_2)dt+\bm f_2\cd dW(t) & \mbox{ in } (0,T)\t G,\\
	\ns\ds w=0 & \mbox{ on } (0,T)\t \G,
\end{cases}	 
\end{equation}
where $ \hat\b, \b $ are given in (\ref{3}), $ f_1\in L^2_\dbF(0,T, L^2(G)) $, $ \bm f_2\in L_\dbF^2(0,T; H^1(G; \dbR^n))  $  and $ \bm f_3\in L^\i_{\dbF}(0,T;$ $L^\i (G; \dbR^n)) $.  Then there exist constants $ \l_0=\l_0(\bm f_3) $ and $ \mu_0=\mu_0(\bm f_3 ) $ such that for any $ \b\in [0,1] $, $ \l\ge \l_0 $  and $ \mu\ge \mu_0 $, we have
\begin{equation}\label{c5}
\ba{ll}
\ds   \frac{1}{2} \dbE\int_0^T\int_{G}  \th^2 \big[ \hat\b^2  |\D w|^2 +\l^2\mu^2 (t+2)^{2\mu-2}  w^2+2\l\mu(t+2)^{\mu-1} |\n w|^2 \big] dxdt\\
\ns \ds + \dbE\int_0^T \int_{G} \Big[\frac{1}{2} \l\mu(t+2)^{\mu-1} \th|\bm f_2|^2+ \frac{1-\b^2}{2} \th^2|\div \bm f_2|^2\Big]dxdt\\
\ns\ds \le 4 \dbE \int_0^T \int_{G} \th^2 f_1^2dxdt+ e^{2\l (T+2)^\mu } \dbE \int_{G}  \big[ \l \mu (T+2)^{\mu-1} w(T,x)^2+\hat \b  |\n w(T,x)|^2 \big]dx,
\ea 
\end{equation}
where the weight function $ \th $ is given in (\ref{4}).
\end{theorem}
\begin{theorem}\label{th4}
Let $ p\in L_\dbF^2 (0,T; H^2(G)) $ be a solution to
\begin{equation}\label{12}
dp -\hat \b\D p  dt=g_1dt- \b \n p\cd dW(t)
\end{equation}
for some $ g_1 \in L_\dbF^2 (0,T; L^2(\dbR^n))$ with the boundary condition $ p=0 $ on $ (0,T)\t \pa G $. The weight function is the same as in (\ref{4}). Then 	for any $ \l\ge 0 $, $ \mu\ge 144(T+2)^2 $ and $ \b\in [0,1]$, it holds that
\begin{eqnarray}\label{c6} 
&& \frac{1}{4}\l\mu^2 \dbE \int_0^T\int_{G}  (t+2)^{\mu-2}\th^2 p^2 dxdt+\sqrt \mu \dbE\int_0^T\int_{G} \th^2 |\n p|^2 dxdt \nonumber \\
&& \le 2 \dbE \int_0^T \int_{G} \th^2g_1^2dxdt+ 2 \l\mu(T+2)^{\mu-1}e^{2\l(T+2)^\mu} \dbE \int_{G} p(T,x)^2dx\\
&&\q +\hat \b e^{2\l2^\mu}\dbE \int_{G}|\n p(0,x)|^2dx+\sqrt \mu e^{2\l 2^{\mu}} \dbE \int_{\dbR^n} p(0,x)^2dx. \nonumber 
\end{eqnarray}
\end{theorem}

The proofs of Theorems \ref{th3} and \ref{th4} are very similar to those for  Theorems \ref{th1} and  \ref{th2}, respectively. Hence, we omit them.

\section{ Proof of Theorems \ref{theorem1} and \ref{theorem2}}


\par {\bf Proof of Theorem \ref{theorem1}. } Let $u=u_1-u_2$, $\rho = \rho_1 -\rho_2$ and $U=U_1-U_2$. Then we have
\begin{eqnarray}\label{b1} 
&& \n \cd \big(\rho_1\n_p B(t,x,\n u_1)\big)-\n\cd \big(\rho_2 \n_p B(t,x,\n u_2)\big) \nonumber\\
&& = \n\cd \big(\rho \n_p B(t,x,\n u_1)\big) + \n \rho_2\cd  \big(\n_p B(t,x,\n u_1)-\n_p B(t,x,\n u_2)\big) \nonumber\\
&& \q +\rho_2 \Big[ \sum_{j=1}^n \big(B_{x_jp_j}(t,x,\n u_1)-B_{x_jp_j} (t,x,\n u_2)\big)\nonumber\\
&& \qq\qq +\sum_{j,k =1}^n  \big( B_{p_jp_k} (t,x,\n u_1)u_{1x_jx_k}-B_{p_jp_k} (t,x,\n u_2)u_{2x_jx_k} \big)\Big]\\
&& =\n\rho\cd \n_p B(t,x,\n u_1)+\rho \(\sum_{j=1}^n B_{p_jx_j} (t,x,\n u_1)+\sum_{j,k=1}^n  B_{p_jp_k} (t,x,\n u_1)u_{1x_{j}x_k}\)\nonumber\\
&& \q +\sum_{j,k =1}^n u_{x_j} \rho_{2x_k} \int_0^1  B_{p_jp_k}\big(t,x, \n u_1+\tau (\n u_2-\n u_1)\big)d\tau \nonumber\\
&& \q +  \sum_{j,k=1}^n u_{x_j} \rho_2 \int_{0}^1 B_{p_jp_kx_k}\big(t,x, \n u_1+\tau (\n u_2-\n u_1)\big)d\tau  +\sum_{j,k=1}^n u_{x_jx_k} \rho_2 B_{p_jp_k} (t,x, \n u_1)\nonumber\\
&& \q +\sum_{j,j',k'=1}^n  u_{x_j} \rho_2 u_{2x_{j'}x_{k'}}  \int_0^1 B_{p_{j'}p_jp_k}\big(t,x, \n u_1+\tau (\n u_2-\n u_1)\big)d \tau \nonumber\\
&& =B_1 \rho+\bm B_2\cd \n \rho+\bm B_3 \cd \n u+\sum_{j,k=1}^n B_{4jk} u_{x_jx_k}, \nonumber
\end{eqnarray}
where
\begin{equation}\label{b2-1*}
\left\{	\ba{ll}
\ds B_1 = \sum_{j=1}^n B_{p_jx_j} (t,x,\n u_1)+\sum_{j,k=1}^n  B_{p_jp_k} (t,x,\n u_1)u_{1x_{j}x_k},\\

\ns\ds \bm B_{2} =\n_p B(t,x,\n u_1),\\
\ns\ds B_{4jk}=\rho_2 B_{p_jp_k} (t,x, \n u_1),\q j,k=1,2,\cdots,n.
\ea \right.	
\end{equation}
and  $ \bm B_3\=(B_{31},\cdots, B_{3j},\cdots, B_{3n}) $, for $ j=1,2,\cdots,n $,
\begin{eqnarray}\label{b2-1}
&& B_{3j}=\sum_{k=1}^n \rho_{2x_k} \int_0^1  B_{p_jp_k}\big(t,x, \n u_1+\tau (\n u_2-\n u_1)\big)d\tau\nonumber\\
&& \qq\q+ \sum_{k=1}^n  \rho_2 \int_{0}^1 B_{p_jp_kx_k}\big(t,x, \n u_1+\tau (\n u_2-\n u_1)\big)d\tau\nonumber\\
&&\qq\q +\sum_{j',k'=1}^n   \rho_2 u_{2x_{j'}x_{k'}}  \int_0^1 B_{p_{j'}p_jp_k}\big(t,x, \n u_1+\tau (\n u_2-\n u_1)\big)d \tau.
\end{eqnarray}
On the other hand, we have
\begin{eqnarray}\label{b1-1}
&& B(t,x, \n u_1)-B(t,x,\n u_2)\nonumber\\
&& = \sum_{j=1}^n u_{x_j} \int_0^1	B_{p_j}\big(t,x, \n u_1+\tau (\n u_2-\n u_1)\big)d\tau \\
&&=\bm B_5\cd \n u,\nonumber
\end{eqnarray}
where
\begin{equation}\label{b1-2}
\bm B_5=\int_0^1 \n_p B\big(t,x, \n u_1+\tau (\n u_2-\n u_1)\big)d\tau.
\end{equation}
Also, we have 
\begin{equation}\label{30}
\ba{ll}
\ds F\(\int_{\dbR^n} K(t, x,y)\rho_1(t,y)dy, \rho_1(t,x)\)-F\( \int_{\dbR^n} K(t, x,y)\rho_2(t,x)dy, \rho_2(t,x)\)\\
\ns\ds =F_1(t,x)\int_{\dbR^n} K(t, x,y)\rho(t,y)dy+F_2(t,x)\rho(t,x),
\ea 
\end{equation}
where
$$
\begin{aligned}
F_1(t,x)=\int_0^1 F_y \Big\{\int_{\dbR^n} K(t, x,y)\big[&\rho_1(t,y)+\tau (\rho_2(t,y)-\rho_1(t,y))\big]dy,\\  &\rho_1(t,x)+\tau \big(\rho_2(t,x)-\rho_1(t,x)\big)\Big\}d\tau
\end{aligned}
$$
and
$$
\begin{aligned}
F_2(t,x)=\int_0^1 F_z \Big\{\int_{\dbR^n} K(t, x,y)\[&\rho_1(t,y)+\tau (\rho_2(t,y)-\rho_1(t,y))\]dy,\\ & \rho_1(t,x)+\tau \big(\rho_2(t,x)-\rho_1(t,x)\big)\Big\}d\tau.
\end{aligned}
$$
Here $ F_y$ and $F_z $ represent the derivatives with respect to the first and second argument of $ F(\cd, \cd) $, respectively. 
By Assumptions \ref{ass1} and \ref{ass2}, there exists a constant $ C_1=C_1(M,F)>0 $ such that  
\begin{eqnarray}\label{31}
&& \Vert  B_1\Vert_{L^\i_\dbF(0,T;L^\i(\dbR^n)} + \sum_{i=2,3,5}\Vert \bm B_i\Vert_{L^\i_\dbF((0,T);L^\i( \dbR^n; \dbR^n)} + \sum_{j,k=1}^n \Vert B_{4jk}\Vert_{L^\i_\dbF(0,T;L^\i \dbR^n)}\nonumber \\
&& + \sum_{i=1,2}\!\Vert F_i\Vert_{L^\i_{\dbF}(0,T;L^\i( \dbR^n))} \le  C_1.
\end{eqnarray}
\par Combining (\ref{b1}), (\ref{30}) with (\ref{a1-11}) and (\ref{3}),  we obtain 
\begin{eqnarray}\label{28}
& d \rho-\hat \b\D \rho dt=&\!\! -\(B_1\rho+\bm B_2\cd \n \rho+\bm B_3\cd \n u+ \sum_{j,k=1}^n B_{4jk} u_{x_jx_k} \)dt \nonumber\\ 
&&\!\!-\b \n \rho \cd dW(t),\qq\qq\qq\q   (t,x)\in (0,T)\t \dbR^n,
\end{eqnarray}
and
\begin{eqnarray}\label{28*}
& du+\hat \b \D udt=&\!\!-\b \div Udt- \bm  B_5\cd \n udt\nonumber \\
&&\!\!-\[F_1(t,x)\int_{\dbR^n} \!  K(t,x,y)\rho(t,y)dy\!+\!F_2(t,x)\rho(t,x)\]dt\\
&&  \!\!+U\cd dW(t), \qq\qq\qq \qq (t,x)\in (0,T)\t \dbR^n.\nonumber
\end{eqnarray}

Applying the Carleman estimate  (\ref{11}) to (\ref{28*}) and recalling (\ref{34*}),  we have
\begin{eqnarray}\label{32-1}
&&\dbE\int_0^T \int_{\dbR^n}   \th^2 \[  \( \sum_{j,k=1}^n u_{x_jx_k}^2\) +\l^2\mu^2 (t+2)^{2\mu-2}  u^2+2\l\mu(t+2)^{\mu-1} |\n u|^2 \] dxdt \nonumber\\
&& \le 4 C_1^2\dbE\int_0^T\int_{\dbR^n} \th^2|\n u|^2dxdt+  4 C_1^2\dbE\int_0^T \int_{\dbR^n} \(\int_{\dbR^n} K(t, x,y)\th(t) \rho(t,y)dy\)^2dxdt\\
&& \q + 4 C_1^2\dbE \!\int_0^T\!\int_{\dbR^n}\! \th^2\rho^2dxdt \!+\! 2 e^{2\l (T+2)^\mu } \dbE \int_0^T\! \int_{\dbR^n}\! \big[\l \mu (T\!+\!2)^{\mu-1} u(T,x)^2+\!\hat \b  |\n u(T,x)|^2 \big]dx.\nonumber
\end{eqnarray}
By Assumption \ref{ass2} and H\"older's inequality, we obtain that
\begin{eqnarray}\label{32*}
&&\int_{\dbR^n} \(\int_{\dbR^n} K(t, x,y)\th(t) \rho(t,y)dy\)^2dx \nonumber\\
&& \le  \int_{\dbR^n} |\th(t)\rho(t,y)|^2 dy  \int_{\dbR^n} \int_{\dbR^n} |K(t,x,y)|^2dydx \\
&& \le M^2 \int_{\dbR^n} |\th(t)\rho(t,y)|^2 dy\nonumber
\end{eqnarray}
This, together with  (\ref{32-1}), implies that
\begin{equation}\label{32}
\ba{ll}
\ds  \dbE\int_0^T \int_{\dbR^n}   \th^2 \big[ \hat\b^2  |\D u|^2 +\l^2\mu^2 (t+2)^{2\mu-2}  u^2+2\l\mu(t+2)^{\mu-1} |\n u|^2 \big] dxdt\\
\ns\ds \le 8C_1^2 M^2\dbE\int_0^T\int_{\dbR^n} \th^2(|\n u|^2+\rho^2)dxdt\\
\ns\ds \q + 2 e^{2\l (T+2)^\mu } \dbE \int_0^T \int_{\dbR^n} \big[ \l \mu (T+2)^{\mu-1} u(T,x)^2+\hat \b  |\n u(T,x)|^2 \big]dx.
\ea 	
\end{equation}

Applying Theorem \ref{th2} to the first equation of (\ref{28}), multiplying each sides of (\ref{20}) by $ \mu^{-1} $, and noting Assumption \ref{ass3},  we get that
\begin{eqnarray}\label{33}
&& \dbE \int_0^T \int_{\dbR^n} \th^2 \[ \frac{1}{4} \l\mu(t+2)^{\mu-2} \rho^2+ \mu^{-\frac{1}{2}} |\n \rho|^2\]dxdt\nonumber \\
&&\le 2\mu^{-1} C_1^2\dbE \int_0^T \int_{\dbR^n} \th^2\(\sum_{j,k=1}^n u_{x_jx_k}^2 +|\n \rho|^2 +|\n u|^2 +\rho^2 \)dxdt\\
&&\q +2\l(T+2)^{\mu-1} e^{2\l(T+2)^{\mu}} \dbE \int_{\dbR^n} \rho(T,x)^2dx+\hat \b \mu^{-1}e^{2\l 2^{\mu}}\dbE \int_{\dbR^n} |\n \rho (0,x) |^2dx\nonumber\\
&&+ e^{2\l 2^{\mu}} \dbE \int_{\dbR^n} p(0,x)^2dx\nonumber.
\end{eqnarray}
By (\ref{32}) and (\ref{33}),  we obtain that
\begin{eqnarray}\label{41}
&&\big(\hat\b^2-2 \mu^{-1}C_1^2\big)\dbE \int_0^T \int_{\dbR^n} \th^2 \( \sum_{j=1}^n u_{x_jx_j}^2\)dxdt+\l^2\mu^2 2^{2\mu-2}\dbE \int_0^T \int_{\dbR^n} \th^2 u^2dxdt \nonumber\\
&&\q +\big(2\l\mu 2^{\mu-1}-8C_1^2M^2-2\mu^{-1}C_1^2 \big)\dbE \int_0^T \int_{\dbR^n} \th^2|\n u|^2dxdt \nonumber\\
&&\q + \(\frac{1}{4}\l\mu 2^{\mu-2}-2\mu^{-1} C_1^2 -8C_1^2M^2\) \dbE \int_0^T \int_{\dbR^n} \th^2\rho^2dxdt\\
&&\q +\big(\mu^{-\frac{1}{2}} -2\mu^{-1}C_1^2\big) \dbE \int_0^T \int_{\dbR^n} \th^2|\n \rho|^2dxdt \nonumber\\
&&\le 2e^{2\l (T+2)^\mu } \dbE \int_0^T \int_{\dbR^n} \big[ \l \mu (T+2)^{\mu-1} u(T,x)^2+\hat \b  |\n u(T,x)|^2 \big]dx \nonumber\\
&&\q +2\l(T+2)^{\mu-1} e^{2\l(T+2)^{\mu}} \dbE  \int_{\dbR^n}  \rho(T,x)^2dx+\hat \b \mu^{-1}e^{2\l 2^{\mu}}\dbE \int_{\dbR^n} |\n \rho(0,x)|^2dx. \nonumber\\
&&\q + e^{2\l 2^{\mu}} \dbE \int_{\dbR^n} p(0,x)^2dx.\nonumber
\end{eqnarray}
By choosing $ \l=1 $, it can be seen that there exists a positive constant $ \mu_0=\mu_0(M, F) $ such that  when $ \mu \geq \mu_0 $,   the terms on the left side of (\ref{41}) are all positive. Consequently, 
there exists a constant $ C=C(T, M,F)>0 $ such that (\ref{27}) holds.
\endpf

\ms

{\bf Proof of Theorem \ref{theorem2}. } Let $u=u_1-u_2$, $\rho = \rho_1 -\rho_2$ and $U=U_1-U_2$. By (\ref{41}),  there exist  positive constants $ \mu_0=\mu_0(T, M, F) $ and $ C=C(T, M, F) $, such that for any $ \l\ge 0 $,     
\begin{equation}\label{a4}
\ba{ll}
\ds \dbE \int_0^T\int_{\dbR^n} \th^2 (\rho^2+|\n \rho|^2 +u^2+|\n u|^2)dxdt\\
\ns\ds \le C e^{2\l (T+2)^{\mu_0}} \dbE \int_{\dbR^n} \big(u(T,x)^2+|\n u(T,x)|^2+ \rho (T,x)^2\big) dx\\
\ns\ds \q + Ce^{2\l 2^{\mu_0}}\dbE \int_{\dbR^n} \big(\rho (0,x)^2+|\n \rho (0,x)|^2\big) dx.
\ea 	
\end{equation}
Recalling that $ \th= e^{\l (t+2)^{\mu_0}} $,   we have
\begin{eqnarray}\label{a5}
&&\Vert \rho\Vert_{L^2_\dbF (\e, T; H^1(\dbR^n))}+\Vert u\Vert_{L^2(\e, T; H^1(\dbR^n))} \nonumber \\
&& \le C e^{\l[(T+2)^{\mu_0}-(2+\e)^{\mu_0}]} \big(\Vert u(T)\Vert_{L^2_{\cF_T}(\O; H^1(\dbR^n))}+\Vert \rho(T)\Vert_{L_{\cF_T}^2(\O; L^2(\dbR^n))}\big)\\
&& \q + Ce^{\l[2^{\mu_0}-(2+\e)^{\mu_0}]}\Vert \rho(0)\Vert_{H^1(\dbR^n)}. \nonumber
\end{eqnarray}
Set  
$$ \l \=\frac{\ln \Vert \rho(0)\Vert_{H^1(\dbR^n)}-\ln \big(\Vert u(T)\Vert_{L^2_{\cF_T}(\O; H^1(\dbR^n))}+\Vert \rho(T)\Vert_{L_{\cF_T}^2(\O; L^2(\dbR^n))}\big)}{(T+2)^{\mu_0}-2^{\mu_0}} $$
in (\ref{a5}). Then we get
\begin{equation*}\label{a6}
\begin{aligned} &
\Vert \rho\Vert_{L^2_\dbF (\e, T; H^1(\dbR^n))}+\Vert u\Vert_{L^2(\e, T; H^1(\dbR^n))} \\ & \leq 	C \Vert \rho(0)\Vert_{H^1(\dbR^n)}^{1-\eta}\big(\Vert u(T)\Vert_{L^2_{\cF_T}(\O; H^1(\dbR^n))}+\Vert \rho(T)\Vert_{L_{\cF_T}^2(\O; L^2(\dbR^n))}\big)^\eta,
\end{aligned}
\end{equation*}
with
\begin{equation*}\label{a7}
\eta =\frac{(2+\e)^{\mu_0}-2^{\mu_0}}{(T+2)^{\mu_0}-2^{\mu_0}}.
\end{equation*}
\endpf

\section{Proof of Theorem \ref{theorem3}}

\par {\bf Proof of Theorem \ref{theorem3}. } Let $ u=u_1-u_2 $, $ \rho=\rho_1-\rho_2 $, $ U=U_1-U_2 $ and $ r=r_{1}-r_{2} $.  Similar to (\ref{28}), the following holds:
\begin{eqnarray}\label{b2}
\left\{\ba{ll}
\ds d \rho-\hat \b\D  \rho dt=  -\(B_1\rho+\bm B_2\cd \n \rho+ \bm B_3\cd \n u+ \sum_{j,k=1}^n B_{4jk}u_{x_jx_j} \)dt  \\
\ns\ds \qq\qq\qq~~-\b \n \rho \cd dW(t) & \mbox{in } Q,\\
\ns\ds du+\hat \b \D udt=\(rR-\b \div U-\bm B_5\cd \n u\)dt&\\
\ns\ds \qq\qq\qq~~ -\[F_1 \int_{G}   K(\cd, \cd,y)\rho(\cd,y)dy+F_2 \rho \]dt +U\cd dW(t) &  \mbox{in } Q,\\
\ns\ds u=\rho=u_{x_1}=\rho_{x_1}=u_{x_1x_1}=\rho_{x_1x_1}=0 & \mbox{on } (0,T)\t \pa G,\\
\ns\ds u(T)=\rho(T)=0, & \mbox{in } G.
\ea \right.
\end{eqnarray}
For arbitrary small $ \e>0 $, we choose $ t_1 $ and $ t_2 $ such that
\begin{equation}\label{b3}
0<t_1<t_2<\e.	
\end{equation}
Let $ \chi\in C^\i(\dbR) $ satisfy that $ 0\le \chi\le 1 $ and that 
\begin{equation}\label{b4}
\chi=\left\{\ba{ll}
\ds 0, &t\le t_1,\\
\ns\ds 1, & t\ge t_2.
\ea \right.	
\end{equation}
Put  $ v= \frac{u}{R} $ and $ V=\frac{V}{R} $.  By elementary calculation, we have that
\begin{eqnarray}\label{c1}
&& dv+\hat\b  \D vdt \nonumber \\
&& = 	\(-\frac{R_t}{R} v+\hat\b \frac{\D R}{R} v+2\hat \b \frac{\n R\cd \n v}{R}\)dt + rdt \nonumber \\
&& \q  + \(-\b\frac{\n R\cd V}{R} -\b \div V- \bm B_5 \cd \n v-\frac{\bm B_5 \cd \n R v}{R} \)dt  \\
&& \q -\frac{1}{R} \( F_1  \int_{\dbR^n} K(\cd,\cd,y)\rho(\cd,y)dy+F_2 \rho \)dt +V \cd dW(t) \nonumber \\
&& = \big(r+f_1 v+ \bm f_2\cd \n v+\bm f_3\cd  V-\b \div V\big)dt \nonumber \\
&& \q -\frac{1}{R} \( F_1 \int_{\dbR^n} K(\cd,\cd, y)\rho(\cd ,y)dy+F_2\rho \)dt +V \cd dW(t),  \nonumber 
\end{eqnarray}
where
\begin{equation}\label{a1-6}
\begin{cases}
\ds f_1\= -\frac{R_t}{R}+\hat\b \frac{\D R}{R}-\frac{\bm B_5\cd \n R}{R},\\
\ns\ds  \bm f_2\=2\hat \b \frac{\n R}{R}-\bm B_5,  \q \bm f_3\=-\b\frac{\n R}{R}. 
\end{cases} 
\end{equation}
Setting $ w\deq\chi v_{x_1} $, $ \cW \deq \chi V_{x_1} $ and $ p\deq\chi \rho_{x_1}  $, and noting that $ r $ is independent of $ x_1 $,   we get from \eqref{b2} and \eqref{c1} that
\begin{eqnarray}\label{b6}
&& d w+\hat \b \D w dt\nonumber\\
&& = \big(\chi_t  v_{x_1}+f_1 w+\bm f_2\cd \n w+ f_{1x_1} \chi v+ \chi  { \bm f}_{2x_1}\cd \n v  \big)dt  + \big(\bm f_3\cd \cW  -\b \div \cW+\chi {\bm f}_{3x_1}\cd V\big)dt \nonumber\\
&& \q -\frac{1}{R} \( F_{1x_1} \int_{G} K(\cd, \cd,y)\chi\rho (\cd,y)dy+F_1 \int_{G} K_{x_1}(\cd,\cd,y)\chi \rho(\cd,y)dy\)dt\\
&&\q -\big(F_{2x_1} \chi\rho  +F_2 p \big)dt +\cW \cd dW (t)  \nonumber
\end{eqnarray}
and that
\begin{equation}\label{c2}
\ba{ll}
\ds dp-\b \D pdt\\
\ns\ds = -\(\chi_t \rho_{x_1}+ g_1 w +\bm g_2 \cd \n w+\sum_{j,k=1}^n g_{3jk}w_{x_jx_k}+\bm g_4 \cd \n p+g_5 p\)dt\\
\ns\ds \q -\chi \[ g_{1x_1} v+ {\bm g_{2x_1}} \cd \n v+\sum_{j,k=1}^n  (\pa_{x_1}g_{3jk})v_{x_jx_k}+ {\bm g_{4x_1}} \cd \n \rho+ g_{5x_1} \rho \]dt+\b \n p\cd dW(t),
\ea 
\end{equation}
where $ \bm g_2= (g_{21},\cdots,g_{2j},\cdots, g_{2n}) $, $ \bm g_4= (g_{41},\cdots,g_{4j},\cdots, g_{4n}) $  and  for $ j,k=1,\cdots,n $,
\begin{equation}\label{a1-8}
\ba{llllll}
\ds g_1= \bm B_3\cd \n R+\sum_{j,k=1}^n B_{4jk}R_{x_jx_k} , &\ds  g_{2j} = R  B_{3j}+2\sum_{k=1}^nB_{4jk} R_{x_jk} , &\ds g_{3jk}= B_{4jk}R, \\
\ns\ds \bm g_4=\bm B_2, &\ds g_5 =B_1.
\ea 
\end{equation}

From the choice of $\chi$ (see \eqref{b4}), we have that
\begin{equation}\label{c7}
\ba{ll}
\ds \int_G (\chi v)^2dx=\int_{\dbR^{n-1}} dx'\int_0^1 \(\int_0^{x_1} w(t,z,x')dz\)^2dx_1	\le \int_G w^2dx.
\ea 
\end{equation} 
Similarly, we can obtain that
\begin{equation}\label{c8}
\begin{cases}
\ds \int_G |\chi V|^2dx=\int_{\dbR^{n-1}} dx'\int_0^1 \(\int_0^{x_1} \cW (t,z,x')dz\)^2dx_1\le \int_G |\cW|^2dx,\\
\ns\ds \int_G (\chi\rho)^2  dx=\int_{\dbR^{n-1}} dx'\int_0^1 \(\int_0^{x_1} p(t,z,x')dz\)^2dx_1\le \int_G p^2dx,
\end{cases}	
\end{equation}
and that for $ j,k=1,\cdots, n $,  
\begin{equation}\label{c8*}
\begin{cases}
\ds \int_G (\chi  v_{x_j})^2dx=\int_{\dbR^{n-1}} dx'\int_0^1 \(\int_0^{x_1} w_{x_j}(t,z,x')dz\)^2dx_1\le \int_G  w_{x_j}^2dx,\\
\ns\ds \int_G (\chi v_{x_jx_k})^2dx=\int_{\dbR^{n-1}} dx'\int_0^1 \(\int_0^{x_1} w_{x_jx_k}(t,z,x')dz\)^2dx_1\le \int_G w_{x_jx_j}^2dx,\\
\ns\ds \int_G (\chi \rho_{x_j})^2dx=\int_{\dbR^{n-1}} dx'\int_0^1 \(\int_0^{x_1} p_{x_j}(t,z,x')dz\)^2dx_1\le \int_G p_{x_j}^2dx.
\end{cases}		
\end{equation}
At last, similar to (\ref{32*}), using H\"older inequality, combining (\ref{c8}), we can get that
\begin{eqnarray}\label{c9*}
	&&\int_G \(\int_{G}K(t,x,y)\chi(t) \rho(t,y)dy\)^2 dx \nonumber\\
	&& \le  \int_{G} |\chi(t)\rho(t,y)|^2 dy \int_{G} \int_{G} K(t,x,y)^2dxdy \\
	&& \le M^2 \int_{G} p(T,y)^2 dy.\nonumber
\end{eqnarray}
On the other hand, recalling  $\Vert K_{x_1}\Vert_{L_\dbF^\i(0,T; L^2(G\t G))}\le M $, we have
\begin{eqnarray}\label{c10}
&& \int_G\(\int_G K_{x_1}(t, x,y)\chi(t) \rho(t,y)dy\)^2dx \nonumber\\
&& \le  \int_{G} |\chi(t)\rho(t,y)|^2 dy \int_{G}\int_{G}  K_{x_1}(t,x,y)^2 dxdy\\
&&\le M^2 \int_G p(T,y)^2dy.\nonumber
\end{eqnarray}

Next, we apply the Carleman estimates (\ref{c5}) and (\ref{c6}) to (\ref{b6}) and (\ref{c2}), respectively. 
From (\ref{c5}), keeping the Assumption \ref{ass3} in mind, and noting  (\ref{c7})--(\ref{c10}), we find that there exist  constants { $ C_2=C_2(T, M,F, R ) $}, $ \l_0=\l_0(T, M, F, R) $ and $ \mu_0=\mu_0(T, M, F, R) $, such that for every $ \l\ge \l_0 $ and $ \mu\ge \mu_0 $, it holds that
\begin{equation}\label{a1-9}
\ba{ll}
\ds \dbE \int_0^T \int_G \th^2 \[ |\D w|^2+\l^2\mu^2 (t+2)^{2\mu-2}w^2+2\l\mu(t+2)^{\mu-1}|\n w|^2\]dxdt\\
\ns\ds \le \dbE \int_0^T \int_G \[8(\chi_t v_{x_1})^2+ C_2 \th^2 p^2 \]dxdt.
\ea 	
\end{equation}
Multiplying each sides of (\ref{c6}) by $ \mu^{-1} $, we get  that there exist constants { $ C_3=C_3(T, M, F, R) $}, $ \l_1 =\l_1(T, M, F, R)$ and $ \mu_1=\mu_1(T, M, F, R) $, such that for each $ \l\ge \l_1 $ and $ \mu\ge \mu_1 $, 
\begin{equation}\label{a1-10}
\ba{ll}
\ds \l\mu \dbE \int_0^T \int_G (t+2)^{\mu-2}\th^2 p^2 dxdt+\sqrt \mu\dbE\int_0^T \int_Q \th^2 |\n p|^2 dxdt\\
\ns\ds \le C_3\dbE\int_0^T \int_G \th^2 \[\mu^{-1} (\chi_t\rho_{x_1})^2 +\mu^{-1} \th^2 \big(w^2+|\n w|^2+|\D w|^2\big)\]dxdt.
\ea 	
\end{equation}
Adding (\ref{a1-9}) and (\ref{a1-10}) together, we obtain that there exist  constants $ \mu_2=\mu_2(T, M, F, R)>0 $, $ \l_2=\l_2(T, M, F, R)>0 $ and {$ C_4=C_4(T,M,F, R)>0 $} such that  for every $ \mu\ge \mu_2 $ and $ \l\ge \l_2 $, 
\begin{equation}\label{b7}
\ba{ll}
\ds \dbE \int_0^T\int_{G} \th^2 \big(p^2+|\n p|^2 +w^2+|\n w|^2\big)dxdt\le C_4\dbE \int_0^T \int_G \th^2 \big(|\chi_t v_{x_1}|^2+|\chi_t p_{x_1}|^2\big) dxdt.
\ea 	
\end{equation}
Recalling  \eqref{b4} for the definition of $ \chi $ and $ \th=e^{\l (t+2)^\mu} $,  we obtain that
\begin{equation*} 
e^{2\l (\e+2)^{\mu}}\dbE\! \int_ \e^T\! \int_G  \big(p^2+|\n p|^2 +w^2+|\n w|^2\big)dxdt\!\le C_4e^{2\l (t_2+2)^\mu} \dbE\! \int_0^T\! \int_G \big(|v_{x_1}|^2+|p_{x_1}|^2\big) dxdt,
\end{equation*}
which yields that
\begin{equation}\label{b8}
e^{2\l [(\e+2)^{\mu}-(t_2+2)^\mu]}\dbE\! \int_ \e^T\! \int_G\!  (p^2+|\n p|^2 +w^2+|\n w|^2)dxdt\le C_4 \dbE \!\int_0^T\! \int_G\! \big(|v_{x_1}|^2+|p_{x_1}|^2 \big)dxdt.
\end{equation}
Letting $ \l $ tend to infinity, recalling the definition of $ \chi $ in (\ref{b4}), and noting that $ p=\chi \rho_{x_1} $ and $ w=\chi v_{x_1} $, we obtain that
\begin{equation}\label{b1-1-2}
	\rho_{x_1}=v_{x_1}=0  \q \mbox{in } [\e,T]\t G,\q \dbP\mbox{-a.s.}
\end{equation}
\par Noting that $ G\=(0,1)\t \dbR^n $,  and $ \rho=u=0 $ on $ [0,T]\t \pa G $, from (\ref{b1-1-2}) and $ v=u/R $,   we get
\begin{equation}\label{b1-1-3}
\rho=u=0  \q \mbox{ in }[\e, T]\t G,\q \dbP\mbox{-a.s.}
\end{equation}
By the second equation in (\ref{b2}), we have
\begin{equation}\label{b1-1-4}
(rR-\b \div U)dt    =U\cd dW(t), \q \mbox{ in }[\e,T]\t G.
\end{equation}
This, together with $ \inf_{x\in \cl G}|R(t,x)|>0 $,   yields $ r=0 $ in $(\e,T)\t G $, $\dbP$-a.s.  Since $ \e >0$ is arbitrary, we get that $ r=0 $ in $ (0,T)\t G $, $\dbP$-a.s.
\endpf


\begin{thebibliography}{10}


\bibitem{MBMF} M. Bardi and M. Fischer, \it On non-uniqueness and uniqueness of solutions in finite-horizon mean field games. \sl ESAIM Control Optim. Calc. Var., \rm {\bf 25}(2019), Paper No. 44.	


\bibitem{ABJFPY} A. Bensoussan, J. Frehse and P. Yam, \it Mean field games and mean field type control theory.   \rm Springer, New York,  2013.

\bibitem{Cardaliaguet2023} P. Cardaliaguet and F. Delarue, \it
Selected topics in mean field games. \sl ICM--International Congress of Mathematicians. \rm Vol. 5.  3660--3703, EMS Press, Berlin, 2023.

\bibitem{PCFDJLPL} P. Cardaliaguet, F. Delarue, J. Lasry and P. Lions, \sl The master equation and the convergence problem in mean field games. \rm Ann. of Math. Stud., 201, Princeton University Press, \rm 2019.





\bibitem{RCFD1} R. Carmona and F. Delarue, \sl Probabilistic theory of mean field games with applications. I. mean field FBSDEs, control, and games. \rm   Springer,  Cham,  \rm 2018.

\bibitem{RCFD2} R. Carmona and F. Delarue, \it Probabilistic theory of mean field games with applications. II. mean field games with common noise and master equations. \rm  Springer, Cham, 2018.





\bibitem{LDWLSO} L. Ding, W. Li, S. Osher and W. Yin, \it A mean field game inverse problem. \sl J. Sci. Comput., \rm {\bf 92}(2022), Paper No. 7, 35 pp.

\bibitem{MDHLGZ} M. Ding, H. Liu and G. Zheng, \it Determining internal topological structures and running cost of mean field games with partial boundary measurement. \sl arXiv 2408.08911v1, \rm 2024.






\bibitem{HCM} M. Huang, R. P. Malham\'e and P. E. Caines, \it Large population stochastic dynamic games: close-loop MCKean-Vlasov system and the nash certainty equivalence principle. \sl Commun. Inf. Syst., \rm  {\bf 6}(2006), 221--252.

\bibitem{Imanuvilov2024} O. Yu. Imanuvilov, H. Liu and M. Yamamoto,  \it Lipschitz stability for determination of states and inverse source problem for the mean field game equations. \sl
Inverse Probl. Imaging, \rm{\bf 18} (2024), 824--859.

\bibitem{PKSASS} P. Kachroo, S. Agarwal and S. Sastry, \it Inverse problem for non-viscous mean field control: example from traffic. \sl IEEE Trans. Autom. Control, \rm  {\bf 61}(2016), 3412--3421.

\bibitem{MVK1} M. V. Klibanov, \it A coefficient inverse problem for the mean field games system, \sl Appl. Math. Optim.,  \rm  {\bf 88}(2023), Paper No. 54pp.


\bibitem{MVK2} M. V. Klibanov, \it The mean field games system: Carleman estimates, Lipschitz stability and uniqueness. \sl J. Inverse Ill-Posed Probl.,  \rm  {\bf 31}(2023), 455--466.

\bibitem{MVKYA} M. V. Klibanov and  Y. Averboukh, \it Lipschitz stability estimate and uniqueness in the retrospective analysis for the mean field games system via two Carleman estimate. \sl SIAM J. Math.  Anal.,  \rm  {\bf 56}(2024), 616--636.



\bibitem{MVKJLHL} M. V. Klibanov, J. Li and H. Liu, \it H\"older stability and uniqueness for the mean field games system via Carleman estimates.  \sl Stud. Appl. Math.,  \rm  {\bf 151}(2023), 1447--1470.

\bibitem{MVKJLHL1} M. V. Klibanov, J. Li and H. Liu, \it On the mean field games system with the lateral Cauchy data via Carleman estimates. \sl J. Inverse Ill-Posed Probl.,  \rm  {\bf 32}(2024), 277--295.

\bibitem{MVKJLZY} M. V. Klibanov, J. Li and Z. Yang, \it Convexification numerical method for the retrospective problem of mean field games. \sl Appl. Math. Optim.,  \rm  {\bf 90}(2024), Paper No. 24 pp. 

\bibitem{VNK} V. N. Kolokoltsov and O. A. Malafeyev, \sl Many agent games in socio-economic systems: corruption, inspection, coalition building, network growth, security. \rm 
Springer, Cham, 2019.



\bibitem{DL} D. Lacker, \it A general characterization of the mean field limit for stochastic differential games. \sl Probab. Theory Relat. Fields, \rm {\bf 165}(2016), 581--648. 

\bibitem{JLPL} J. M. Larsy and P. L. Lions, \it  Mean field games. \sl Jpn, J. Math.,  \rm {\bf 2}(2007), 229--260.

\bibitem{Liao2024} Z. Liao and Q. L\"u, \it Stability estimate for an inverse stochastic parabolic problem of determining unknown time-varying boundary. \sl
Inverse Problems, \rm {\bf 40} (2024), Paper No. 045032.

\bibitem{HLCMSZ} H. Liu, C. Mou and S. Zhang, \it Inverse Problems for mean field games. \sl Inverse Problems, \rm  {\bf 39}(2023), Paper No. 085003.

\bibitem{QL} Q. L\"u, \it Carleman estimate for stochastic parabolic equations and inverse stochastic parabolic problems. \sl Inverse Problems, \rm {\bf 28}(2012), 045008.

\bibitem{QLXZ} Q. L\"u, X. Zhang, \it Inverse problems for stochastic partial differential equations: some progresses and open problems. \sl Numer. Algebra Control Optim., \rm 14(2024), 227--272.

\bibitem{XZQL} 	Q. L\"u and X. Zhang, \sl Mathematical control theory for stochastic partial differential equations. \rm Springer, Cham,  \rm 2022.






\bibitem{KRNSKW} K. Ren, N. Soedjak and K. Wang, \it Unique determination of cost functions in a multipopulation mean field game model. \sl airXiv: 2312.01622, \rm 2023.

\bibitem{KRNSKWHZ} K. Ren, N. Soedjak, K. Wang and  H. Zhai, \it Reconstructing a state-independent cost function in a mean field game model. \sl arXiv: 2402.09297, \rm 2024. 

\bibitem{STXZ} S. Tang and X. Zhang, \it Null controllability for forward and backward stochastic parabolic equations. \sl SIAM J. Control Optim.,  \rm  {\bf 48}(2009), 2196--2216.

\bibitem{Wu2020} B. Wu, Q. Chen and Z. Wang, \it Carleman estimates for a stochastic degenerate parabolic equation and applications to null controllability and an inverse random source problem.
\sl
Inverse Problems \rm {\bf 36} (2020),  075014.

\bibitem{Wu2022} B. Wu and J. Liu, \it On the stability of recovering two sources and initial status in a stochastic hyperbolic-parabolic system. \sl
Inverse Problems \rm {\bf 38} (2022),  Paper No. 025010.

\bibitem{GY} G. Yuan, \it Conditional stability in determination of initial data for stochastic parabolic equations. \sl Inverse Problems, \rm 33(2017), Paper No. 035014.













\end{thebibliography}
\end{document}